\documentclass[twoside,11pt,amsfonts]{article} 

\newcommand{\abst}[1]{\medskip {\footnotesize #1  } \medskip }

\newcommand{\sect}[1]{\bigskip{\bf #1}\medskip}

\newcommand{\definition}[1]{\medskip \bf #1 \rm}

\newcommand{\bth}[1]{\medskip \bf #1 \it}
\newcommand{\eth}{\rm \smallskip}

\newcommand{\remark}[1]{\medskip {\bf Remark\, #1} \rm}

\newcommand{\example}[1]{\medskip {\bf Example\, #1} \rm}

\newcommand{\proof}{\rm P\,r\,o\,o\,f. }
\newcommand{\proofend}{\hfill \hbox{\vrule width 5pt height 5pt depth 0pt}}

\newfont{\sml}{cmsl10}
\newfont{\sss}{msam10}
\newfont{\Bbb}{msbm10} 
\newfont{\goth}{eufm10} 
\newfont{\cysc}{wncysc10} 

\newcommand{\R}{\mbox{\Bbb R}}

\newcommand{\head}[5]
{
\thispagestyle{empty}

\begin{flushright}
{\footnotesize
This is a preprint of a paper whose final and definite form 
has been published in:\\ Math. Balkanica {\bf 26} (2012), no~1-2, 191--202.}
\end{flushright}

\markboth{\hfill \rm #5} {\rm #4 \hfill}

\vspace*{0.8cm}

\begin{flushleft}
{\Large\bf #1}  

\bigskip
\smallskip

{\Large\it #2} 

\bigskip
{\small \it #3} 

\end{flushleft}
\bigskip
}

\usepackage{amsmath}
\usepackage{cite}

\begin{document}

\head
 {Calculus of Variations with Classical \\ [4pt] and Fractional Derivatives}
 {Tatiana Odzijewicz $^1$, Delfim F. M. Torres $^2$}
 {Presented at 6$^{th}$ International Conference ``TMSF' 2011''}
 {Calculus of Variations with Classical \dots}
 {T. Odzijewicz, D.F.M. Torres}


\abst{
We give a proper fractional extension of the classical calculus of variations.
Necessary optimality conditions of Euler--Lagrange type for variational problems
containing both classical and fractional derivatives are proved.
The fundamental problem of the calculus of variations
with mixed integer and fractional order derivatives
as well as isoperimetric problems are considered.

\medskip

{\sl MSC 2010}: 49K05, 26A33

{\sl Key Words}: variational analysis, optimality,
Riemann--Liouville fractional operators,
fractional differentiation, isoperimetric problems
}

\smallskip


\sect{1. Introduction}

One of the classical problems of mathematics consists in finding a
closed plane curve of a given length that encloses the greatest
area: the \emph{isoperimetric problem}. The legend says that the
first person who solved the isoperimetric problem was Dido, the
Queen of Carthage, who was offered as much land as she could
surround with the skin of a bull. Dido's problem is nowadays part of
the \emph{calculus of variations} \cite{MR0160139,Bruce_van_Brunt}.

Fractional calculus is a generalization of (integer) differential calculus,
allowing to define derivatives (and integrals) of real or complex order
\cite{kilbas,miller,I_Podlubny}. The first application of fractional calculus
belongs to Niels Henrik Abel (1802--1829) and goes back to 1823 \cite{Abel}.
Abel applied the fractional calculus to the solution of an integral equation
which arises in the formulation of the \emph{tautochrone problem}.
This problem, sometimes also called the \emph{isochrone problem},
is that of finding the shape of a frictionless wire lying in a vertical
plane such that the time of a bead placed on the wire slides
to the lowest point of the wire in the \emph{same time} regardless
of where the bead is placed. The cycloid is the isochrone
as well as the \emph{brachistochrone} curve:
it gives the \emph{shortest time} of slide and marks
the born of the \emph{calculus of variations}.

The study of fractional problems of the calculus of variations and
respective Euler--Lagrange type equations is a subject of current
strong research due to its many applications in science and
engineering, including mechanics, chemistry, biology, economics, and
control theory \cite{TM}. In 1996--1997 Riewe obtained a version of
the Euler--Lagrange equations for fractional variational problems
combining the conservative and nonconservative cases
\cite{CD:Riewe:1996,CD:Riewe:1997}. Since then, numerous works on
the fractional calculus of variations, fractional optimal control
and its applications have been written---see, \textrm{e.g.},
\cite{MR2356715,MyID:210,Almeida1,Almeida3,MyID:172,Atanackovic,Baleanu,El-Nabulsi:Torres,Frederico:Torres1,Frederico:Torres2,Klimek,Malinowska}
and references therein. For the study of fractional isoperimetric
problems, see \cite{Isoperimetric}.

In the pioneering paper \cite{OmPrakashAgrawal}, and others that followed,
the fractional necessary optimality conditions are proved under the
hypothesis that admissible functions $y$ have continuous left
and right fractional derivatives on the closed interval $[a,b]$.
By considering that the admissible functions $y$ have continuous left fractional derivatives
on the whole interval, then necessarily $y(a)=0$; by considering that the admissible
functions $y$ have continuous right fractional derivatives, then necessarily $y(b)=0$.
This fact has been independently remarked, in different contexts, at least in
\cite{Isoperimetric,Almeida3,Atanackovic,Jelicic}.

In our work we want to be able to consider arbitrarily given boundary conditions
$y(a)=y_a$ and $y(b)=y_b$ (and isoperimetric constraints).
For that we consider variational functionals with integrands involving
not only a fractional derivative of order $\alpha \in (0,1)$
of the unknown function $y$, but also the classical derivative $y'$.
More precisely, we consider dependence of the integrands on the independent
variable $t$, unknown function $y$, and $y'+ k \, _{a}\textsl{D}_t^\alpha y$
with $k$ a real parameter. As a consequence, one gets a proper extension
of the classical calculus of variations, in the sense that the classical
theory is recovered with the particular situation $k = 0$. We remark
that this is not the case with all the previous literature on
the fractional variational calculus, where the classical theory
is not included as a particular case and only as a limit,
when $\alpha \rightarrow 1$.

The text is organized as follows. In Section~2
we briefly recall the necessary definitions and properties
of the fractional calculus in the sense of Riemann--Liouville.
Our results are stated, proved, and illustrated through an example,
in Section~3. We end with Section~4 of conclusion.


\sect{2. Preliminaries}
\label{sec:prelim}

In this section some basic definitions and properties of fractional
calculus are given. For more on the subject we refer the reader to
the books \cite{kilbas,miller,I_Podlubny} and historical survey
\cite{TM}.

\definition{Definition 1.} (Left and right Riemann--Liouville derivatives)
Let $f$ be a function defined on $[a,b]$. The operator
$_{a}\textsl{D}_t^\alpha$, \vskip -8pt
\begin{equation*}
_{a}\textsl{D}_t^\alpha f(t)
=\frac{1}{\Gamma(n-\alpha)}D^n\int_a^t(t-\tau)^{n-\alpha-1}f(\tau)d\tau\, ,
\end{equation*}
is called the left Riemann--Liouville fractional derivative of order
$\alpha$, and the operator $_{t}\textsl{D}_b^\alpha$, \vskip -10pt
\begin{equation*}
_{t}\textsl{D}_b^\alpha f(t)
=\frac{-1}{\Gamma(n-\alpha)}D^n\int_t^b(\tau-t)^{n-\alpha-1}f(\tau)d\tau\, ,
\end{equation*}
is called the right Riemann--Liouville fractional derivative of order $\alpha$,
where $\alpha \in \R^+$ is the order of the derivatives
and the integer number $n$ is such that $n-1\leq\alpha<n$.

\definition{Definition 2.} (Mittag--Leffler function)
Let $\alpha, \beta>0$. The Mittag--Leffler function is defined by
\vskip -10pt
\begin{equation*}
E_{\alpha,\beta}(z)
=\sum_{k=0}^\infty\frac{z^k}{\Gamma(\alpha k+\beta)}\, .
\end{equation*}

\bth{Theorem 3.} (Integration by parts) If $f,g$ and the fractional
derivatives $_{a}\textsl{D}_t^\alpha g$ and $_{t}\textsl{D}_b^\alpha
f$ are continuous at every point $t\in[a,b]$, then \vskip -10pt
\begin{equation}
\label{eq:ip}
\int_a^b f(t)_{a}\textsl{D}_t^\alpha g(t)dt
=\int_a^b g(t)_{t}\textsl{D}_b^\alpha f(t)dt
\end{equation}
\vskip -3pt \noindent
 for any $0<\alpha<1$. \eth

\remark{4.}
\label{rem:ir:bc}
If $f(a)\neq 0$, then $\left._{a}\textsl{D}_t^\alpha f(t)\right|_{t=a}=\infty$.
Similarly, if $f(b)\neq 0$, then $\left._{t}\textsl{D}_b^\alpha f(t)\right|_{t=b}=\infty$.
Thus, if $f$ possesses continuous left and right Riemann--Liouville fractional derivatives
at every point $t\in[a,b]$, then $f(a)=f(b)=0$. This explains why the usual term
$\left. f(t) g(t) \right|_{a}^{b}$ does not appear on the right-hand side of \eqref{eq:ip}.

\medskip


\sect{3. Main results} \label{sec:mr}

Following \cite{Jelicic}, we prove optimality conditions
of Euler--Lagrange type for variational problems containing
classical and fractional derivatives simultaneously.
In Section~3.1 the fundamental variational problem
is considered, while in Section~3.2 we study
the isoperimetric problem. Our results cover
fractional variational problems subject
to arbitrarily given boundary conditions.
This is in contrast with
\cite{OmPrakashAgrawal,MR2356049,MR2356715,MR2519151},
where the necessary optimality conditions are valid
for appropriate zero valued boundary conditions
(\textrm{cf.} Remark~4).
For a discussion on this matter
see \cite{Almeida3,Atanackovic,Jelicic}.


\sect{3.1. The Euler--Lagrange equation}
\label{sub:sec:bp}

Let $0<\alpha<1$. Consider the following problem:
find a function $y \in C^1[a,b]$ for which the functional
\begin{equation}
\label{eq:funct}
\mathcal{J}(y)=\int_a^b F\left(t,y(t),y'(t)+ k \, _{a}\textsl{D}_t^\alpha y(t)\right) dt
\end{equation}
subject to given boundary conditions
\begin{equation}
\label{eq:bc}
y(a)=y_a, \quad y(b)=y_b,
\end{equation}
has an extremum. We assume $k$ is a fixed real number,
$F\in C^2([a,b]\times\R^2;\R)$,
and $\partial_3F$ (the partial derivative of $F(\cdot,\cdot,\cdot)$
with respect to its third argument) has a continuous right Riemann--Liouville
fractional derivative of order $\alpha$.

\definition{Definition 5.}
A function $y \in C^1[a,b]$ that satisfies
the given boundary conditions \eqref{eq:bc} is said to be
\emph{admissible} for problem \eqref{eq:funct}--\eqref{eq:bc}.

\medskip

For simplicity of notation we introduce the operator $[\cdot]_k^\alpha$ defined by
\begin{equation*}
[y]_k^\alpha(t) = \left(t,y(t),y'(t)+ k \, _{a}\textsl{D}_t^\alpha y(t)\right) \, .
\end{equation*}
With this notation we can write \eqref{eq:funct} simply as
\begin{equation*}
\mathcal{J}(y) = \int_a^b F[y]_k^\alpha(t) dt \, .
\end{equation*}

\bth{Theorem 6.} (The fractional Euler--Lagrange equation)
If $y$ is an extremizer (minimizer or maximizer)
of problem \eqref{eq:funct}--\eqref{eq:bc},
then $y$ satisfies the Euler--Lagrange equation
\begin{equation}
\label{eq:4}
\partial_2 F[y]_k^\alpha(t) - \frac{d}{dt}\partial_3 F[y]_k^\alpha(t)
+ k \, _{t}\textsl{D}_b^\alpha \partial_3 F[y]_k^\alpha(t) = 0
\end{equation}
for all $t\in[a,b]$.
\eth

\proof
Suppose that $y$ is a solution of \eqref{eq:funct}--\eqref{eq:bc}.
Note that admissible functions $\hat{y}$ can be written
in the form $\hat{y}(t)=y(t)+\epsilon\eta(t)$, where
$\eta \in C^1[a,b]$, $\eta(a)=\eta(b)=0$,
and $\epsilon\in\R$. Let
$$
J(\epsilon)
=\int_a^b F\left(t,y(t)+\epsilon\eta(t),\frac{d}{dt}\left(y(t)
+\epsilon\eta(t)\right)+ k _{a}\textsl{D}_t^\alpha \left(y(t)
+\epsilon\eta(t)\right)\right)dt.
$$
Since $_{a}\textsl{D}_t^\alpha$ is a linear operator, we know that
$$
_{a}\textsl{D}_t^\alpha\left(y(t)+\epsilon\eta(t)\right)
={_{a}\textsl{D}_t^\alpha} y(t)+\epsilon_{a}\textsl{D}_t^\alpha\eta(t).
$$
On the other hand,
\begin{equation}
\label{eq:1}
\begin{split}
\left.\frac{dJ}{d\epsilon}\right|_{\epsilon=0}
&=\left.\int_a^b \frac{d}{d\epsilon}
F[\hat{y}]_k^\alpha(t) dt\right|_{\epsilon=0}\\
&=\int_a^b\Bigl(\partial_2 F[y]_k^\alpha(t) \cdot \eta(t)
+\partial_3 F[y]_k^\alpha(t) \frac{d\eta(t)}{dt}
+k\partial_3 F[y]_k^\alpha(t) _{a}\textsl{D}_t^\alpha\eta(t)\Bigr)dt.
\end{split}
\end{equation}
Using integration by parts we get
\begin{equation}
\label{eq:2}
\int_a^b\partial_3F\frac{d\eta}{dt}dt=\left.\partial_3
F\eta\right|_a^b - \int_a^b \left(\eta\frac{d}{dt}\partial_3F\right)dt
\end{equation}
and
\begin{equation}
\label{eq:3}
\int_a^b\partial_3F _{a}\textsl{D}_t^\alpha\eta dt
=\int_a^b \eta \, {_{t}}\textsl{D}_b^\alpha \partial_3F dt.
\end{equation}
Substituting \eqref{eq:2} and \eqref{eq:3} into \eqref{eq:1},
and having in mind that $\eta(a)=\eta(b)=0$, it follows that
\begin{equation*}
\left.\frac{dJ}{d\epsilon}\right|_{\epsilon=0}
=\int_a^b \eta(t) \Bigl(\partial_2 F[y]_k^\alpha(t)
-\frac{d}{dt}\partial_3 F[y]_k^\alpha(t)
+ k\, _{t}\textsl{D}_b^\alpha \partial_3 F[y]_k^\alpha(t)\Bigr)dt.
\end{equation*}
A necessary optimality condition is given by
$\left.\frac{dJ}{d\epsilon}\right|_{\epsilon=0}=0$. Hence,
\begin{equation}
\label{eq:bef:afl}
\int_a^b \eta(t)\Bigl(\partial_2 F[y]_k^\alpha(t)
-\frac{d}{dt}\partial_3 F[y]_k^\alpha(t)
+k \, _{t}\textsl{D}_b^\alpha\partial_3 F[y]_k^\alpha(t)\Bigr)dt=0.
\end{equation}
We obtain equality \eqref{eq:4} by applying the fundamental lemma
of the calculus of variations to \eqref{eq:bef:afl}.
\proofend

\medskip

\example{7.}
Note that for $k=0$ our necessary optimality condition \eqref{eq:4}
reduces to the classical Euler--Lagrange equation
\cite{MR0160139,Bruce_van_Brunt}.

\medskip


\sect{3.2. The fractional isoperimetric problem}
\label{sub:sec:fip}

As before, let $0<\alpha<1$.
We now consider the problem
of extremizing a functional
\begin{equation}
\label{eq:funct3}
\mathcal{J}(y)=\int_a^b F\left(t,y(t),y'(t)
+ k\, _{a}\textsl{D}_t^\alpha y(t)\right)dt
\end{equation}
in the class $y \in C^1[a,b]$
when subject to given boundary conditions
\begin{equation}
\label{eq:funct4}
y(a)=y_a\, , \quad y(b)=y_b,
\end{equation}
and an isoperimetric constraint
\begin{equation}
\label{eq:funct2}
\mathcal{I}(y)
=\int_a^b G(t,y(t),y'(t)+ k\, _{a}\textsl{D}_t^\alpha y(t))dt=\xi \, .
\end{equation}
We assume that $k$ and $\xi$ are fixed real numbers,
$F, G\in C^2([a,b]\times \R^2;\R)$,
and $\partial_3F$ and $\partial_3G$ have continuous
right Riemann--Liouville fractional derivatives of order $\alpha$.

\definition{Definition 8.}
A function $y \in C^1[a,b]$ that satisfies
the given boundary conditions \eqref{eq:funct4} and isoperimetric
constraint \eqref{eq:funct2} is said to be \emph{admissible}
for problem \eqref{eq:funct3}--\eqref{eq:funct2}.

\definition{Definition 9.}
An admissible function $y$ is an \emph{extremal}
for $\mathcal{I}$ if it satisfies
the fractional Euler--Lagrange equation
\begin{equation*}
\partial_2 G[y]_k^\alpha(t)
-\frac{d}{dt}\partial_3 G[y]_k^\alpha(t)
+ k\, _{t}\textsl{D}_b^\alpha \partial_3 G[y]_k^\alpha(t) = 0
\end{equation*}
for all $t\in[a,b]$.

\medskip

The next theorem gives a necessary optimality condition
for the fractional isoperimetric problem
\eqref{eq:funct3}--\eqref{eq:funct2}.

\bth{Theorem 10.}
\label{thm:noc:ip}
Let $y$ be an extremizer to the functional \eqref{eq:funct3}
subject to the boundary conditions \eqref{eq:funct4} and
the isoperimetric constraint \eqref{eq:funct2}.
If $y$ is not an extremal for $\mathcal{I}$,
then there exists a constant $\lambda$ such that
\begin{equation}
\label{eq:30}
\partial_2 H[y]_k^\alpha(t)
-\frac{d}{dt}\partial_3 H[y]_k^\alpha(t)
+k \, {_{t}\textsl{D}_b^\alpha} \partial_3 H[y]_k^\alpha(t)=0
\end{equation}
for all $t\in[a,b]$, where $H(t,y,v)=F(t,y,v)-\lambda G(t,y,v)$.
\eth

\proof
We introduce the two parameter family
\begin{equation}
\label{eq:25}
\hat{y}=y+\epsilon_1\eta_1+\epsilon_2\eta_2,
\end{equation}
in which $\eta_1$ and $\eta_2$ are such that $\eta_1,\eta_2\in C^1[a,b]$
and they have continuous left and right fractional derivatives.
We also require that
$$
\eta_1(a)=\eta_1(b)=0=\eta_2(a)=\eta_2(b).
$$
First we need to show that in the family \eqref{eq:25}
there are curves such that $\hat{y}$ satisfies \eqref{eq:funct2}.
Substituting $y$ by $\hat{y}$ in \eqref{eq:funct2},
$\mathcal{I}(\hat{y})$ becomes a function of two parameters
$\epsilon_1,\epsilon_2$. Let
$$
\hat{I}(\epsilon_1,\epsilon_2)
=\int_a^b G(t,\hat{y},\hat{y}'+ k _{a}\textsl{D}_t^\alpha \hat{y})dt-\xi.
$$
Then, $\hat{I}(0,0)=0$ and
$$
\left.\frac{\partial\hat{I}}{\partial\epsilon_2}\right|_{(0,0)}
=\int_a^b\eta_2\left(\partial_2G-\frac{d}{dt}\partial_3G
+k _{t}\textsl{D}_b^\alpha\partial_3G\right)dt.
$$
Since $y$ is not an extremal for $\mathcal{I}$, by the fundamental
lemma of the calculus of variations there is a function $\eta_2$ such that
$$
\left.\frac{\partial\hat{I}}{\partial\epsilon_2}\right|_{(0,0)}\neq 0.
$$
By the implicit function theorem, there exists a function
$\epsilon_2(\cdot)$ defined in a neighborhood of zero, such that
$\hat{I}(\epsilon_1,\epsilon_2(\epsilon_1))=0$.
Let $\hat{J}(\epsilon_1,\epsilon_2)=\mathcal{J}(\hat{y})$.
Then, by the Lagrange multiplier rule, there exists a real $\lambda$ such that
$$
\nabla(\hat{J}(0,0)-\lambda\hat{I}(0,0))=\textbf{0}.
$$
Because
$$
\left.\frac{\partial\hat{J}}{\partial\epsilon_1}\right|_{(0,0)}
=\int_a^b\eta_1\left(\partial_2F-\frac{d}{dt}\partial_3F
+k_{t}\textsl{D}_b^\alpha\partial_3F\right)dt
$$
and
$$
\left.\frac{\partial\hat{I}}{\partial\epsilon_1}\right|_{(0,0)}
=\int_a^b\eta_1\left(\partial_2G-\frac{d}{dt}\partial_3G
+k_{t}\textsl{D}_b^\alpha\partial_3G\right)dt,
$$
one has
\begin{equation*}
\int_a^b\eta_1\Biggl[\left(\partial_2F-\frac{d}{dt}\partial_3F
+k_{t}\textsl{D}_b^\alpha\partial_3F\right)
-\lambda\left(\partial_2G-\frac{d}{dt}\partial_3G
+k_{t}\textsl{D}_b^\alpha\partial_3G\right)\Biggr]dt=0.
\end{equation*}
Since $\eta_1$ is an arbitrary function, \eqref{eq:30} follows
from the fundamental lemma of the calculus of variations.
\proofend


\pagebreak

\sect{3.3. An example}
\label{sub:sec:ex}

Let $\alpha\in\left(0,1\right)$ and $k, \xi \in\R$.
Consider the following fractional isoperimetric problem:
\begin{equation}
\label{eq:ex}
\begin{gathered}
\mathcal{J}(y)=\int_0^1\left(y' + k\, {_{0}\textsl{D}_t^\alpha} y\right)^2dt \longrightarrow \min\\
\mathcal{I}(y)=\int_0^1\left(y'+ k \, {_{0}\textsl{D}_t^\alpha} y\right)dt = \xi\\
y(0)=0 \, , \
y(1)=\int_0^1 E_{1-\alpha,1}\left(-k\left(1-\tau\right)^{1-\alpha}\right) \xi d\tau.
\end{gathered}
\end{equation}
In this case the augmented Lagrangian $H$ of Theorem~10 is given by
$H(t,y,v) =v^2 -\lambda v$. One can easily check that
\begin{equation}
\label{eq:y:ex}
y(t)=\int_0^t E_{1-\alpha,1}\left(-k\left(t-\tau\right)^{1-\alpha}\right) \xi d\tau
\end{equation}
\begin{itemize}
\item is not an extremal for $\mathcal{I}$;
\item satisfies $y'+ k \,{_{0}\textsl{D}_t^\alpha} y= \xi$
(see, \textrm{e.g.}, \cite[p.~297, Theorem~5.5]{kilbas}).
\end{itemize}
Moreover, \eqref{eq:y:ex} satisfies \eqref{eq:30} for $\lambda=2\xi$, \textrm{i.e.},
\begin{equation*}
-\frac{d}{dt}\left(2\left(y'+ k \,{_{0}\textsl{D}_t^\alpha} y\right) -2\xi\right)\\
+k \, {_{t}\textsl{D}_1^\alpha}\left(2\left(y'
+ k \, {_{0}\textsl{D}_t^\alpha} y\right) -2\xi\right)=0.
\end{equation*}
We conclude that \eqref{eq:y:ex} is the extremal for problem \eqref{eq:ex}.

\example{11.}
Choose $k=0$. In this case the isoperimetric
constraint is trivially satisfied,
\eqref{eq:ex} is reduced to the
classical problem of the calculus of variations
\begin{equation}
\label{eq:k0}
\begin{gathered}
\mathcal{J}(y)=\int_0^1 (y'(t))^2 dt \longrightarrow \min\\
y(0)=0 \, , \quad y(1)= \xi,
\end{gathered}
\end{equation}
and our general extremal \eqref{eq:y:ex} simplifies
to the well-known minimizer $y(t)=\xi t$ of \eqref{eq:k0}.

\example{12.}
When $\alpha \rightarrow 1$ the isoperimetric
constraint is redundant with the boundary conditions,
and the fractional problem \eqref{eq:ex} simplifies
to the classical variational problem
\begin{equation}
\label{eq:ex:alpha1}
\begin{gathered}
\mathcal{J}(y)=(k+1)^2 \int_0^1 y'(t)^2 dt \longrightarrow \min\\
y(0)=0 \, , \quad y(1)= \frac{\xi}{k+1}.
\end{gathered}
\end{equation}
Our fractional extremal \eqref{eq:y:ex} gives
$y(t)=\frac{\xi}{k+1}t$, which is exactly
the minimizer of \eqref{eq:ex:alpha1}.

\example{13.}
Choose $k = \xi = 1$. When $\alpha \rightarrow 0$
one gets from \eqref{eq:ex} the classical isoperimetric problem
\begin{equation}
\label{eq:ex:alpha0}
\begin{gathered}
\mathcal{J}(y)=\int_0^1\left(y'(t) + y(t)\right)^2 dt \longrightarrow \min\\
\mathcal{I}(y)=\int_0^1 y(t) dt = \frac{1}{e}\\
y(0)=0 \, , \quad y(1)= 1-\frac{1}{\mathrm{e}}.
\end{gathered}
\end{equation}
Extremal \eqref{eq:y:ex} is then reduced to the classical
extremal $y(t)=1 - \mathrm{e}^{-t}$ of \eqref{eq:ex:alpha0}.

\example{14.}
Choose $k=1$ and $\alpha=\frac{1}{2}$.
Then \eqref{eq:ex} gives the following
fractional isoperimetric problem:
\begin{equation}
\label{eq:ex:alpha=1/2}
\begin{gathered}
\mathcal{J}(y)=\int_0^1\left(y' +  {_{0}\textsl{D}_t^\frac{1}{2}} y\right)^2 dt\longrightarrow \min\\
\mathcal{I}(y)=\int_0^1\left(y'+ {_{0}\textsl{D}_t^\frac{1}{2}} y\right)dt=\xi \\
y(0) = 0 \, , \quad y(1) = -\xi\left(1- \mathrm{erfc}(1)+\frac{2}{\sqrt{\pi}}\right),
\end{gathered}
\end{equation}
where $\mathrm{erfc}$ is the complementary error function.
The extremal \eqref{eq:y:ex} for the particular fractional problem \eqref{eq:ex:alpha=1/2} is
$$
y(t)=-\xi\left(1- \mathrm{e}^t \mathrm{erfc}(\sqrt{t})+\frac{2\sqrt{t}}{\sqrt{\pi}}\right)\, .
$$


\sect{4. Conclusion}
\label{sec:conc}

Fractional variational calculus provides a very useful framework to
deal with nonlocal dynamics in Mechanics and Physics
\cite{R:A:D:10,Baleanu20101111}. It has received considerable
interest in recent years, with several researchers applying this
field to develop fractional classical and quantum mechanics
\cite{cresson-ho,MR2549615,MalinowskaTorres:Hahn}. Motivated by the
results and insights of \cite{Isoperimetric,Almeida1,Jelicic}, in
this paper we generalize previous fractional Euler--Lagrange
equations by proving optimality conditions for fractional problems
of the calculus of variations where the highest derivative in the
Lagrangian is of integer order. This approach avoids the
difficulties with the given boundary conditions when in presence of
the Riemann--Liouville derivatives \cite{Jelicic}. For the case with
the Caputo fractional derivatives (\cite{I_Podlubny}) we refer the
reader to \cite{withTatiana:Basia}.

We focus our attention to problems subject
to integral constraints (fractional isoperimetric problems),
which have recently found a broad class of important applications
\cite{isoJMAA,Viktor,Curtis}. For $k = 0$ our results are
reduced to the classical ones \cite{Bruce_van_Brunt}.
This is in contrast with the standard approach to fractional variational calculus,
where the integer-order case is obtained only in the limit.


\medskip

{\bf Acknowledgements.}\, This work was first announced at the IFAC
Workshop on Fractional Derivatives and Applications (IFAC FDA'2010),
held in University of Extremadura, Badajoz, Spain, October 18-20,
2010; then subsequently at conference ``TMSF' 2011''. It was
supported by {\it FEDER} funds through {\it COMPETE}
--- Operational Programme Factors of Competitiveness (``Programa
Operacional Factores de Competitividade'') and by Portuguese funds
through the {\it Center for Research and Development in Mathematics
and Applications} (CIDMA), University of Aveiro, and the Portuguese
Foundation for Science and Technology (``FCT --- Funda\c{c}\~{a}o
para a Ci\^{e}ncia e a Tecnologia''), within project
PEst-C/MAT/UI4106/2011 with COMPETE number
FCOMP-01-0124-FEDER-022690. Odzijewicz was also supported by FCT
through the Ph.D. fellowship SFRH/BD/33865/2009.



\bigskip

\it

\noindent
$^{1,2}$ Center for Research and Development in Mathematics and Applications \\
Department of Mathematics, University of Aveiro \\
3810-193 Aveiro, PORTUGAL \\[4pt]
e-mail: $^1$ tatianao@ua.pt , $^2$ delfim@ua.pt \hfill Received:
October 13, 2011



{\small

\begin{thebibliography}{99}

\bibitem{Abel}
N.H. Abel,
{\it Euvres Completes de Niels Henrik Abel}.
Christiana: Imprimerie de Grondahl and Son, New York and London,
Johnson Reprint Corporation (1965).

\vspace*{-8pt}

\bibitem{OmPrakashAgrawal}
O.P. Agrawal,
Formulation of {E}uler-{L}agrange equations for fractional variational problems.
{\it J. Math. Anal. Appl.} {\bf 272}, No~1 (2002), 368--379.

\vspace*{-8pt}

\bibitem{MR2356049}
O.P. Agrawal,
A general finite element formulation for fractional variational problems.
{\it J. Math. Anal. Appl.} {\bf 337}, No~1 (2008), 1--12.

\vspace*{-8pt}

\bibitem{MR2356715}
O.P. Agrawal, D. Baleanu,
A {H}amiltonian formulation and a direct numerical scheme
for fractional optimal control problems.
{\it J. Vib. Control} {\bf 13}, No~9-10 (2007), 1269--1281.

\vspace*{-8pt}

\bibitem{Isoperimetric}
R. Almeida, R.A.C. Ferreira, D.F.M. Torres,
Isoperimetric problems of the calculus of variations with fractional derivatives.
{\it Acta Math. Sci. Ser. B Engl. Ed.} {\bf 32}, No~2 (2012), 619--630.
{\tt arXiv:1105.2078}

\vspace*{-8pt}

\bibitem{R:A:D:10}
R. Almeida, A.B. Malinowska, D.F.M. Torres,
A fractional calculus of variations for multiple integrals
with application to vibrating string.
{\it J. Math. Phys.} {\bf 51}, No~3 (2010), 033503, 12pp.
{\tt arXiv:1001.2722}

\vspace*{-8pt}

\bibitem{MyID:210}
R. Almeida, S. Pooseh, D.F.M. Torres,
Fractional variational problems depending on indefinite integrals.
{\it Nonlinear Anal.} {\bf 75}, No~3 (2012), 1009--1025.
{\tt arXiv:1102.3360}

\vspace*{-8pt}

\bibitem{Almeida1}
R. Almeida, D.F.M. Torres,
Calculus of variations with fractional derivatives and fractional integrals.
{\it Appl. Math. Lett.} {\bf 22}, No~12 (2009), 1816--1820.
{\tt arXiv:0907.1024}

\vspace*{-8pt}

\bibitem{isoJMAA}
R. Almeida, D.F.M. Torres,
H\"olderian variational problems subject to integral constraints.
{\it J. Math. Anal. Appl.} {\bf 359}, No~2 (2009), 674--681.
{\tt arXiv:0807.3076}

\vspace*{-8pt}

\bibitem{cresson-ho}
R. Almeida, D.F.M. Torres,
Generalized Euler-Lagrange equations
for variational problems with scale derivatives.
{\it Lett. Math. Phys.} {\bf 92}, No~3 (2010), 221--229.
{\tt arXiv:1003.3133}

\vspace*{-8pt}

\bibitem{Almeida3}
R. Almeida, D.F.M. Torres,
Leitmann's direct method for fractional optimization problems.
{\it Appl. Math. Comput.} {\bf 217}, No~3 (2010), 956--962.
{\tt arXiv:1003.3088}

\vspace*{-8pt}

\bibitem{MyID:172}
R. Almeida, D.F.M. Torres,
Fractional variational calculus for nondifferentiable functions.
{\it Comput. Math. Appl.} {\bf 61}, No~10 (2011), 3097--3104.
{\tt arXiv:1103.5406}

\vspace*{-8pt}

\bibitem{Atanackovic}
T.M. Atanackovi{\'c}, S. Konjik, S. Pilipovi{\'c},
Variational problems with fractional derivatives: {E}uler-{L}agrange equations.
{\it J. Phys. A} {\bf 41}, No~9 (2008), 095201, 12pp.
{\tt arXiv:1101.2961}

\vspace*{-8pt}

\bibitem{Baleanu}
D. Baleanu,
New applications of fractional variational principles.
{\it Rep. Math. Phys.} {\bf 61}, No~2 (2008), 199--206.

\vspace*{-8pt}

\bibitem{MR2519151}
D. Baleanu, O. Defterli, O.P. Agrawal,
A central difference numerical scheme for fractional optimal control problems.
{\it J. Vib. Control} {\bf 15}, No~4 (2009), 583--597.

\vspace*{-8pt}

\bibitem{Baleanu20101111}
D. Baleanu, J.I. Trujillo,
A new method of finding the fractional Euler-Lagrange and Hamilton
equations within Caputo fractional derivatives.
{\it Commun. Nonlinear Sci. Numer. Simul.} {\bf 15}, No~5 (2010), 1111--1115.

\vspace*{-8pt}

\bibitem{Viktor}
V. Bl{\aa}sj\"{o},
The isoperimetric problem.
{\it Amer. Math. Monthly} {\bf 112}, No~6 (2005), 526--566.

\vspace*{-8pt}

\bibitem{MR2549615}
J. Cresson, G.S.F. Frederico, D.F.M. Torres,
Constants of motion for non-differentiable quantum variational problems.
{\it Topol. Methods Nonlinear Anal.} {\bf 33}, No~2 (2009), 217--231.
{\tt arXiv:0805.0720}

\vspace*{-8pt}

\bibitem{Curtis}
J.P. Curtis,
Complementary extremum principles for isoperimetric optimization problems.
{\it Optim. Eng.} {\bf 5}, No~4 (2004), 417--430.

\vspace*{-8pt}

\bibitem{El-Nabulsi:Torres}
R.A. El-Nabulsi, D.F.M. Torres,
Necessary optimality conditions for fractional action-like integrals of variational
calculus with {R}iemann-{L}iouville derivatives of order {$(\alpha,\beta)$}.
{\it Math. Methods Appl. Sci.} {\bf 30}, No~15 (2007), 1931--1939.
{\tt arXiv:math-ph/0702099}

\vspace*{-8pt}

\bibitem{Frederico:Torres1}
G.S.F. Frederico, D.F.M. Torres,
A formulation of {N}oether's theorem for fractional problems
of the calculus of variations.
{\it J. Math. Anal. Appl.} {\bf 334}, No~2 (2007), 834--846.
{\tt arXiv:math/0701187}

\vspace*{-8pt}

\bibitem{Frederico:Torres2}
G.S.F. Frederico, D.F.M. Torres,
Fractional conservation laws in optimal control theory.
{\it Nonlinear Dynam.} {\bf 53}, No~3 (2008), 215--222.
{\tt arXiv:0711.0609}

\vspace*{-8pt}

\bibitem{MR0160139}
I.M. Gelfand, S.V. Fomin,
{\it Calculus of Variations}.
Revised English edition translated and edited by Richard A.
Silverman, Prentice-Hall Inc., Englewood Cliffs, N.J. (1963).

\vspace*{-8pt}

\bibitem{Jelicic}
Z.D. Jelicic, N. Petrovacki,
Optimality conditions and a solution scheme for fractional optimal control problems.
{\it Struct. Multidiscip. Optim.} {\bf 38}, No~6 (2009), 571--581.

\vspace*{-8pt}

\bibitem{kilbas}
A.A. Kilbas, H.M. Srivastava, J.J. Trujillo,
{\it Theory and Applications of Fractional Differential Equations}.
Elsevier Science B.V., Amsterdam (2006).

\vspace*{-8pt}

\bibitem{Klimek}
M. Klimek,
Lagrangean and {H}amiltonian fractional sequential mechanics.
{\it Czechoslovak J. Phys.} {\bf 52}, No~11 (2002), 1247--1253.

\vspace*{-8pt}

\bibitem{TM}
J.T. Machado, V. Kiryakova, F. Mainardi,
Recent history of fractional calculus.
\emph{Commun. Nonlinear Sci. Numer. Simul.}
{\bf 16}, No~3 (2011), 1140--1153.

\vspace*{-8pt}

\bibitem{Malinowska}
A.B. Malinowska, D.F.M. Torres,
Generalized natural boundary conditions for fractional variational
problems in terms of the Caputo derivative.
{\it Comput. Math. Appl.} {\bf 59}, No~9 (2010), 3110--3116.
{\tt arXiv:1002.3790}

\vspace*{-8pt}

\bibitem{MalinowskaTorres:Hahn}
A.B. Malinowska, D.F.M. Torres,
The Hahn quantum variational calculus.
{\it J. Optim. Theory Appl.} {\bf 147}, No~3 (2010), 419--442.
{\tt arXiv:1006.3765}

\vspace*{-8pt}

\bibitem{miller}
K.S. Miller, B. Ross,
{\it An Introduction to the Fractional Calculus and Fractional Differential Equations}.
A Wiley-Interscience Publication, John Wiley \& Sons Inc., New York (1993).

\vspace*{-8pt}

\bibitem{withTatiana:Basia}
T. Odzijewicz, A.B. Malinowska, D.F.M. Torres,
Fractional variational calculus with classical and combined Caputo derivatives.
\emph{Nonlinear Anal.} {\bf 75}, No~3 (2012), 1507--1515.
{\tt arXiv:1101.2932}

\vspace*{-8pt}

\bibitem{I_Podlubny}
I. Podlubny,
{\it Fractional Differential Equations}.
Academic Press Inc., San Diego, CA (1999).

\vspace*{-8pt}

\bibitem{CD:Riewe:1996}
F. Riewe,
Nonconservative {L}agrangian and {H}amiltonian mechanics.
{\it Phys. Rev. E (3)} {\bf 53}, No~2 (1996), 1890--1899.

\vspace*{-8pt}

\bibitem{CD:Riewe:1997}
F. Riewe,
Mechanics with fractional derivatives.
{\it Phys. Rev. E (3)} {\bf 55}, No~3, part B (1997), 3581--3592.

\vspace*{-8pt}

\bibitem{Bruce_van_Brunt}
B. van Brunt,
{\it The Calculus of Variations}.
Springer-Verlag, New York (2004).

\end{thebibliography}
}
\end{document}